\documentclass[11pt]{article}

\usepackage{url,epsf,amsmath,amssymb,amsfonts,graphicx,tikz}
\usepackage{color}
\usepackage{hyperref}

\newtheorem{theorem}{\bf Theorem}[section]
\newtheorem{corollary}[theorem]{\bf Corollary}

\newtheorem{proposition}[theorem]{\bf Proposition}

\newtheorem{problem}[theorem]{\bf Problem}

\newcommand{\proof}{\noindent{\bf Proof.\ }}
\newcommand{\qed}{\hfill $\Box$ \bigskip}
\newcommand{\iztok}{\color{black}}
\newcommand{\s}{\color{black}}
\newcommand{\igy}{\color{black}}

\begin{document}

\title{{\s Graphs that are simultaneously efficient open domination and efficient closed domination graphs}}

\author{
Sandi Klav\v zar $^{a,b,c}$
\and
Iztok Peterin $^{c,d}$
\and
Ismael G. Yero $^e$
}

\date{}

\maketitle

\begin{center}
$^a$ Faculty of Mathematics and Physics, University of Ljubljana, Slovenia \\
\medskip

$^b$ Faculty of Natural Sciences and Mathematics, University of Maribor, Slovenia \\
\medskip

$^{c}$ Institute of Mathematics, Physics and Mechanics, Ljubljana, Slovenia \\
\medskip

$^{d}$ Faculty of Electrical Engineering and Computer Science, University of Maribor, Slovenia \\
\medskip

$^{e}$ Departamento de Matem\'aticas, EPS, Universidad de C\'adiz, Algeciras, Spain
\end{center}

\begin{abstract}
A graph is an efficient open (resp.\ closed) domination graph if there
exists a subset of vertices whose open (resp.\ closed) neighborhoods
partition its vertex set. Graphs that are efficient open as well as efficient
closed (shortly EOCD graphs) are investigated. {\s The structure of EOCD graphs with respect to their efficient open and efficient closed dominating sets is explained.} It is shown that the decision problem {\igy regarding} whether a graph is an EOCD graph is an NP-complete problem.  A recursive description that constructs all EOCD trees is given and EOCD graphs {\s are characterized among the Sierpi\'nski graphs}.
\end{abstract}

\noindent
{\bf E-mails}: sandi.klavzar@fmf.uni-lj.si, iztok.peterin@um.si, ismael.gonzalez@uca.es
\medskip

\medskip

\medskip\noindent
{\bf AMS Subj. Class.}: 05C69, 05C05, 68R10, 94B60

\section{Introduction}

The \emph{domination number}, $\gamma (G)$, of a graph $G$ is an important classical
graph invariant with many applications. It is defined as the minimum cardinality of a
subset of vertices $S$, called \emph{dominating set}, with the property that each vertex from $V(G)-S$ has a neighbor in $S$. A dominating set $S$ of cardinality $\gamma (G)$ is called a $\gamma (G)$-\emph{set}. The union of closed neighborhoods centered in vertices of {\s a dominating set} covers the entire vertex set. A classical question for a cover of a set is: when {\igy does this cover form} a partition? A
graph $G$ is called an \emph{efficient closed domination graph}, or {\em ECD
graph} for short, if there exist{\s s} a set $P$, $P\subseteq V(G)$, such that the
closed neighborhoods centered in vertices of $P$ partition $V(G)$. Such a set
$P$ is called a\emph{\ perfect code} of $G$. More general,
a set $P$ is an $r$-\emph{perfect code} of $G$ if the $r$-balls
centered in vertices of $P$ partition $V(G)$.

The study of perfect codes in graphs was initiated by Biggs~\cite{biggs-1973} and presents a generalization of the problem of the existence of (classical) error-correcting codes. The initial research focused on distance regular and related classes of graphs, while later the investigation was extended to general graphs, cf.~\cite{kratochvil-1986}. To determine whether a given graph has a 1-perfect code is an NP-complete problem~\cite{bange-1988} and remains NP-complete on $k$-regular graphs ($k\geq 4$)~\cite{kratochvil-1988}, on planar graphs of maximum degree {\s $3$}~\cite{fellows-1991, kratochvil-1988}, as well as on bipartite and chordal graphs~\cite{smart-1995}. On the positive side, the existence of a 1-perfect code can be decided in polynomial time on trees~\cite{fellows-1991}, interval graphs~\cite{kratochvil-1995}, and circular-arc graphs~\cite{klostermeyer-2000}.

In the last period, the study of perfect codes in graphs was primarily focused on their existence and construction in some central families of graphs. {\s Many researches were} done on standard graph products and product-like graphs~\cite{abay-2009, hrastnik-2015, klavzar-2006, mollard-2011, taylor-2009, zerovnik-2008}. Among other classes of graphs on which perfect codes were investigated we mention Sierpi\'nski graphs~\cite{cull-1999, klavzar-2002}, cubic vertex-transitive graphs~\cite{knor-2012}, circulant graphs~\cite{deng-2014}, {\s twisted tori~\cite{jha-2014}, dual cubes~\cite{jha-2015}}, and {AT}-free and dually chordal graphs~\cite{brandstadt-2015}.

A graph invariant closely related to the domination number is the \emph{total domination number} $\gamma _{t}(G)$~\cite{henning-2013}. It is defined as the minimum cardinality of a subset of vertices $D$, called \emph{total dominating set}, such that each vertex from $V(G)$ has a
neighbor in $D$. A total dominating set $D$ of cardinality $\gamma _{t}(G)$
is called a $\gamma _{t}(G)$-\emph{set}. If we switch to neighborhoods, the union of open
neighborhoods centered in vertices of {\s a total dominating set} covers the entire vertex set {\s and} again one can pose the question: when does this cover form a partition? A graph $G$ is called an \emph{efficient open domination graph}, or an {\em EOD
graph} for short, if there exists a set $D$, $D\subseteq V(G)$, such that open
neighborhoods centered in vertices of $D$ partition $V(G)$. Such a set $D$ is called an \emph{EOD set}. Note that two different
vertices of an EOD set are either adjacent or at distance at least {\s $3$}.

The problem of deciding whether a graph $G$ is an EOD graph is $NP$-complete~\cite{GaScSl,McRae}. For various properties of EOD graphs see~\cite{GaSc}, a recursive characterization of EOD trees is given in~\cite{GaScSl}. EOD graphs that are also Cayley graphs were studied in~\cite{Tham}, while EOD grid graphs were investigated in \cite{CoHeKe, Dej, KlGo}. EOD direct product graphs were characterized in~\cite{AbHamTay}, for other standard graph products (lexicographic, strong, disjunctive and Cartesian) see~\cite{KuPeYe1}. {\s Domination-type problems studied on graph products are usually most difficult on the Cartesian product, recall the famous Vizing's conjecture \cite{BrDoKl}. It is hence not} surprising that EOD graphs studied on product graphs seems to be the most difficult on the Cartesian product. For some very recent results in this direction see~\cite{KPRT}.

In this paper we study the graphs that are ECD and EOD at the same time {\s and call them}{\em efficient open closed domination graphs}, {\em EOCD graphs} for short. In the rest of the paper we shall use the term \textit{ECD set} instead of 1-perfect code to make the notation consistent.

We proceed as follows. In the rest of this section additional definitions {\s are given and a basic result recalled}. In the next section we {\s show ho{\igy w} to construct an ECD graph from and EOD graph and vice versa, and} consider the structure of EOCD graphs from the viewpoint of the relation{\igy ship} between selected EOD sets and selected ECD sets. In two extremal cases we find that for the corresponding EOCD graphs $G$ we have $\gamma_t(G) = \gamma(G)$ and  $\gamma_t(G) = 2\gamma(G)$, respectively. In Section~\ref{sec:complexity} we prove that the decision problem {\igy regarding} whether a graph is an EOCD graph is an NP-complete problem. On the other hand, in one {\s of the above extremal cases}, EOCD graphs can be recognized in polynomial time. Then, in Section~\ref{sec:trees}, we give a recursive description of EOCD trees, while in the final section EOCD graphs {\s are characterized among the Sierpi\'nski graphs}.

We will use the notation $[n] = \{1,\ldots ,n\}$ and $[n]_0 = \{0,\ldots ,n-1\}$. Throughout the article we consider only finite, simple graphs. If $S$ is a subset of vertices of a graph, then $\left\langle S\right\rangle$ denotes the subgraph induced by $S$.  A \emph{matching} of a graph is an independent set of its edges. For the later use we next state the following basic result. {\s Its first assertion has been independently discovered several times, cf~\cite[Theorem 4.2]{haynes-1998}, while for the second fact see~\cite{KuPeYe1}.}

\begin{proposition}
\label{prp:EOD-ECD}
{\s Let $G$ be a graph.}

(i) If $P$ is an ECD set of $G$, then $|P| = \gamma (G)$.

(ii) If $D$ is an EOD set of $G$, then $|D| = \gamma _{t}(G)$.
\end{proposition}

\section{On the structure of EOCD graphs}
\label{sec:structure}

{\s In this section we first show that each EOD graph naturally yields an {\igy E}CD graph and that each ECD graph can be modified to an} {\iztok EOD} {\s graph. Then we consider the structure of EOCD graph{\igy s} with respect to the relation{\igy ship} between their EOD in ECD sets.}

If $D$ is an EOD set of a graph $G$, then $D$ induces a matching $M$. Note that an edge from $M$ lies in no triangle, hence its contracting produces no parallel edges. Now, let $G'$ be the graph obtained from $G$ by contraction all the edges from $M$. Then $G'$ is an ECD graph with an ECD set consisting of the vertices obtained by the contraction of $M$.

Conversely, let $G'$ be an ECD graph with an ECD set $P$. For every vertex $v\in P$ partition {\s the set of its neighbors} arbitrarily into sets $A$ and $B$. (If the degree of $v$ is $1$, then necessarily one of these sets is empty.) Let $G$ be the graph obtained from $G'$ by replacing  every vertex $v\in P$ by two adjacent vertices $v_A$ and $v_B$, and
adding edges $uv_A$ for every $u\in A$ and edges $uv_B$ for every $u\in B$. Then $G$ is an EOD graph with an EOD set $\{v_A,v_B:\ v\in P\}$.

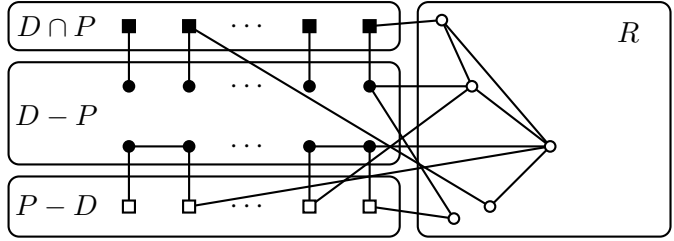
\begin{figure}[htb]
\begin{center}
\begin{tikzpicture}[scale=.8,style=thick,x=1cm,y=1cm]
\def\vr{2.5pt} 
\path (2,0.7) coordinate (v1);
\path (3,0.7) coordinate (v2);
\path (5,0.7) coordinate (v3);
\path (6,0.7) coordinate (v4);
\path (2,1.7) coordinate (u1);
\path (3,1.7) coordinate (u2);
\path (5,1.7) coordinate (u3);
\path (6,1.7) coordinate (u4);
\path (2,2.7) coordinate (w1);
\path (3,2.7) coordinate (w2);
\path (5,2.7) coordinate (w3);
\path (6,2.7) coordinate (w4);
\path (2,3.7) coordinate (z1);
\path (3,3.7) coordinate (z2);
\path (5,3.7) coordinate (z3);
\path (6,3.7) coordinate (z4);
\path (7.2,3.8) coordinate (c);
\path (7.4,0.5) coordinate (c1);
\path (7.7,2.7) coordinate (c2);
\path (9,1.7) coordinate (c3);
\path (8,0.7) coordinate (c4);

%
\draw (z1) -- (w1);
\draw (z2) -- (w2);
\draw (z3) -- (w3);
\draw (z4) -- (w4);
\draw (v1) -- (u1);
\draw (v2) -- (u2);
\draw (v3) -- (u3);
\draw (v4) -- (u4);
\draw (u1) -- (u2);
\draw (u3) -- (u4);
\draw (c) -- (z4);
\draw (c1) -- (w4);
\draw (c1) -- (v4);
\draw (c2) -- (w4);
\draw (c2) -- (v3);
\draw (c3) -- (u4);
\draw (c3) -- (v2);
\draw (c4) -- (z2);
\draw (c2) -- (c3);
\draw (c3) -- (c4);
\draw (c2) -- (c);
\draw (c3) -- (c);


\draw (c) [fill=white] circle (\vr);
\draw (c1) [fill=white] circle (\vr);
\draw (c2) [fill=white] circle (\vr);
\draw (c3) [fill=white] circle (\vr);
\draw (c4) [fill=white] circle (\vr);
\draw (u1) [fill=black] circle (\vr);
\draw (u2) [fill=black] circle (\vr);
\draw (u3) [fill=black] circle (\vr);
\draw (u4) [fill=black] circle (\vr);
\draw (w1) [fill=black] circle (\vr);
\draw (w2) [fill=black] circle (\vr);
\draw (w3) [fill=black] circle (\vr);
\draw (w4) [fill=black] circle (\vr);
\draw [fill=black] (1.9,3.8) rectangle (2.1,3.6);
\draw [fill=black] (2.9,3.8) rectangle (3.1,3.6);
\draw [fill=black] (4.9,3.8) rectangle (5.1,3.6);
\draw [fill=black] (5.9,3.8) rectangle (6.1,3.6);
\draw [fill=white] (1.9,0.8) rectangle (2.1,0.6);
\draw [fill=white] (2.9,0.8) rectangle (3.1,0.6);
\draw [fill=white] (4.9,0.8) rectangle (5.1,0.6);
\draw [fill=white] (5.9,0.8) rectangle (6.1,0.6);
\draw [rounded corners] (0,1.2) rectangle (6.5,0.2);
\draw [rounded corners] (0,3.1) rectangle (6.5,1.4);
\draw [rounded corners] (0,4.1) rectangle (6.5,3.3);
\draw [rounded corners] (6.8,4.1) rectangle (11,0.2);

\draw (4,0.7) node {$\cdots$};
\draw (4,1.7) node {$\cdots$};
\draw (4,2.7) node {$\cdots$};
\draw (4,3.7) node {$\cdots$};
\draw (0.8,0.7) node {$P-D$};
\draw (0.8,2.2) node {$D-P$};
\draw (0.8,3.7) node {$D\cap P$};
\draw (10.3,3.6) node {$R$};

\end{tikzpicture}
\end{center}
\caption{Structure of an EOCD graph.}
\label{fig:structure}
\end{figure}

Let $G$ be an EOCD graph with an EOD set $D$ and an ECD set $P$. Then $V(G)$ can be partitioned into sets $D\cap P$, $D-P$, $P-D$, and $R=V(G)-(D\cup P)$, see Fig.~\ref{fig:structure}. Clearly, some of these sets may be empty. From the definitions of ECD and EOD sets we infer the following properties.

\begin{itemize}
\item A vertex from $D\cap P$ (a black squared vertex in Fig.~\ref{fig:structure}) can have an arbitrary number of neighbors in $R$, has a unique neighbor in $D-P$, and has no neighbors in $P-D$.
\item A vertex from $P-D$ (a white squared vertex in Fig.~\ref{fig:structure}) can have  an arbitrary number of neighbors in $R$, a unique neighbor in $D-P$, and no neighbors in $D\cap P$.
\item A vertex from $D-P$ (a black vertex in Fig.~\ref{fig:structure}) can have an arbitrary number of neighbors in $R$ and, either {\igy a unique neighbor in $P-D$ and a unique neighbor in $D-P$, or a unique neighbor in $D\cap P$}.
\item A vertex from $R$ (a white vertex in Fig.~\ref{fig:structure}) can have an arbitrary number of neighbors in $R$ and either {\igy a unique neighbor in $P-D$ and a unique neighbor in $D-P$, or a unique neighbor in $D\cap P$}.
\item Vertices from $D\cap P$ together with their unique neighbors from $D-P$ induce a matching.
\item Vertices from $P-D$ together with their unique neighbors in $D-P$ induce $k$ copies of $P_{4}$, where $2k=|P-D|$.
\end{itemize}

The described structure is visible in Fig.~\ref{fig:structure}. The same notation will be used later in Fig.~\ref{operation123}.

The described structure above yields two extreme cases: either $D\cap P=\emptyset$ or $P-D=\emptyset$. Clearly, the structure of any EOCD graph depends on the chosen ECD set $P$ and EOD set $D$. That is, different pairs of sets $P,D$ could produce different configurations. In this sense, if there exists an ECD set $P$ and an EOD set $D$ in $G$, such that $D\cap P=\emptyset$, then we say that $G$ is an EOCD graph {\em with empty $D\cap P$}, and if $P-D=\emptyset$, then we say that  $G$ is an EOCD graph {\em with empty $P-D$}. We observe that $D-P$ is always {\s non-empty} for every ECD set $P$ and every EOD set $D$ of any EOCD graph. Moreover, if $R=\emptyset$, then $G$ is formed only {\s from the} disjoint union of copies of $K_2$ {\s and} copies of $P_4$.

The following two propositions follow directly from the above mentioned structure. The first result characterizes the EOCD graphs with empty $D\cap P$.

\begin{proposition}
\label{prp:DcapP:empty}
A graph $G$ is an EOCD graph with empty $D\cap P$ if and only if there exists $A\subseteq V(G)$, such that
$\left\langle A\right\rangle =kP_4$, where every vertex from $V(G)-A$ is adjacent to exactly one vertex
of degree $1$ in $\left\langle A\right\rangle$ and one vertex of degree $2$ in $\left\langle A\right\rangle$.
\end{proposition}

The second result characterizes the EOCD graphs with empty $P-D$.

\begin{proposition}
\label{extreme1}
\label{prp:P-D:empty}
A graph $G$ is an EOCD graph with empty $P-D$ if and only if there exists $D{\igy \subseteq V(G)}$ that induces a matching $M$, where every edge of $M$ contains at least one vertex of degree $1$ in $G$
$($this vertex is from $D-P)$ and every vertex  from $V(G)-D$ is adjacent to exactly one vertex in $M$ which is in $P$.
\end{proposition}

We end this section with a connection between $\gamma (G)$ and $\gamma_t(G)$ for EOCD graphs with empty $D\cap P$ or empty $P-D$, respectively. Both results follow from the described structure of EOCD graphs, and by applying Proposition~\ref{prp:EOD-ECD}.

\begin{proposition}
\label{prp:gamma_t-small}
If $G$ is an EOCD graph with empty $D\cap P$, then $\gamma_t(G)=\gamma (G)$.
\end{proposition}

\begin{proposition}
\label{prp:gamma}
If $G$ is an EOCD graph with empty $P-D$, then $\gamma_t(G)=2\gamma (G)$.
\end{proposition}

{\s Recall that for any graph $G$ (without isolated vertices) $\gamma(G)\le \gamma_t(G)\le 2\gamma(G)$ holds. The above} two propositions are of interest {\s particularly} because it is an open problem to characterize {\igy the} graphs $G$ with $\gamma_t(G) = 2\gamma(G)$, as well as {\igy the} graphs $G$ with $\gamma_t(G) = \gamma(G)$, cf.~\cite[p.~36]{henning-2013}. {\igy In this direction, the} trees $T$ for which $\gamma_t(T) = \gamma(T)$ holds were characterized in~\cite[Theorem 6]{dorfling-2006} as the trees obtained from a disjoint union of $P_4$s by means of certain four operations. {\igy Moreover, a} characterization of trees $T$ for which $\gamma_t(T) = 2\gamma(T)$ holds was obtained in~\cite{henning-2001}. For {\s these two} results see also~\cite[Sections 4.6 and 4.7]{henning-2013}.

\section{Complexity results}
\label{sec:complexity}

In this section we deal with the problem of deciding whether a given graph
contains an EOD set and an ECD set ({\sc EOCD Problem}
for short), that is, the following problem.

\begin{center}
\fbox{\parbox{0.85\linewidth}{\noindent
{\sc EOCD Problem}\\[.8ex]
\begin{tabular*}{.93\textwidth}{rl}
{\em Input:} & A simple graph $G$.\\
{\em Question:} & Is $G$ an EOCD graph?
\end{tabular*}
}}
\end{center}

In order to study {\s this problem}, we shall make a reduction from the
one-in-three 3-SAT problem, which is known to be NP-complete~\cite{garey} {\s and reads as follows}.

\begin{center}
\fbox{\parbox{0.85\linewidth}{\noindent
{\sc One-In-Three 3-SAT}\\[.8ex]
\begin{tabular*}{.93\textwidth}{rl}
{\em Input:} & A Boolean formula $\mathcal{F}$ on $n$ variables and $m$
clauses.\\
{\em Question:} & Is there a satisfying truth assignment for the $n$
variables, \\
& such that each clause has exactly one true literal?
\end{tabular*}
}}
\end{center}

Next we present the main result of this section, which is in part inspired by the proof of the NP-completeness of the problem of deciding
whether a graph contains an {\iztok EOD} set given in~\cite{GaScSl}.

\begin{theorem}
\label{thm:NP-complete}
The {\sc EOCD Problem} is NP-complete.
\end{theorem}

\begin{figure}[ht!]
\centering
\begin{tikzpicture}[scale=0.8]
\node [draw, shape=circle,black,fill=white,scale=0.6] (t1) at  (-2,0.5) {};
\draw (-2.2,0.5) node[left] {$u_i$};
\node [draw, shape=circle,black,fill=white,scale=0.6] (t2) at  (2,0.5) {};
\draw (2.2,0.5) node[right] {$\overline{u_i}$};
\node [draw, shape=circle,black,fill=white,scale=0.6] (t3) at  (0,1.5) {};
\draw (0,1.3) node[below] {$t_{i1}$};

\node [draw, shape=circle,black,fill=white,scale=0.6] (p1) at  (0,2.7) {};
\draw (-0.1,2.7) node[left] {$t_{i2}$};
\node [draw, shape=circle,black,fill=white,scale=0.6] (p2) at  (-2,3.8) {};
\draw (-2,3.6) node[below] {$t_{i3}$};
\node [draw, shape=circle,black,fill=white,scale=0.6] (p3) at  (0,3.8) {};
\draw (-0.35,3.9) node[above] {$q_i$};
\node [draw, shape=circle,black,fill=white,scale=0.6] (p4) at  (2,3.8) {};
\draw (2,3.6) node[below] {$t_{i4}$};

\node [draw, shape=circle,black,fill=white,scale=0.6] (c1) at  (0,5) {};
\draw (0,5.15) node[above] {$c_{i1}$};
\node [draw, shape=circle,black,fill=white,scale=0.6] (c2) at  (-1.5,5) {};
\draw (-1.65,5) node[left] {$c_{i7}$};
\node [draw, shape=circle,black,fill=white,scale=0.6] (c3) at  (1.5,5) {};
\draw (1.65,5) node[right] {$c_{i2}$};
\node [draw, shape=circle,black,fill=white,scale=0.6] (c4) at  (-2,6) {};
\draw (-2.15,6) node[left] {$c_{i6}$};
\node [draw, shape=circle,black,fill=white,scale=0.6] (c5) at  (2,6) {};
\draw (2.15,6) node[right] {$c_{i3}$};
\node [draw, shape=circle,black,fill=white,scale=0.6] (c6) at  (-0.7,7) {};
\draw (-0.85,7.2) node[left] {$c_{i5}$};
\node [draw, shape=circle,black,fill=white,scale=0.6] (c7) at  (0.7,7) {};
\draw (0.85,7.2) node[right] {$c_{i4}$};

\node [draw, shape=circle,black,fill=white,scale=0.6] (f1) at  (-4,8) {};
\draw (-4.15,8) node[left] {$v_{i1}$};
\node [draw, shape=circle,black,fill=white,scale=0.6] (f2) at  (4,8) {};
\draw (4.15,8) node[right] {$v_{i2}$};

\node [draw, shape=circle,black,fill=white,scale=0.6] (g1) at  (-3,9) {};
\draw (-2.8,8.8) node[below] {$w_{i1}$};
\node [draw, shape=circle,black,fill=white,scale=0.6] (g2) at  (3,9) {};
\draw (2.8,8.8) node[below] {$w_{i6}$};
\node [draw, shape=circle,black,fill=white,scale=0.6] (g3) at  (-1.5,9.5) {};
\draw (-1.5,9.3) node[below] {$w_{i3}$};
\node [draw, shape=circle,black,fill=white,scale=0.6] (g4) at  (0,9.5) {};
\draw (0,9.3) node[below] {$w_{i4}$};
\node [draw, shape=circle,black,fill=white,scale=0.6] (g5) at  (1.5,9.5) {};
\draw (1.5,9.3) node[below] {$w_{i5}$};
\node [draw, shape=circle,black,fill=white,scale=0.6] (g6) at  (-3,10) {};
\draw (-3.2,10) node[left] {$w_{i2}$};
\node [draw, shape=circle,black,fill=white,scale=0.6] (g7) at  (3,10) {};
\draw (3.2,10) node[right] {$w_{i7}$};

\draw(t1)--(t2)--(t3)--(t1)--(f1)--(f2)--(t2);
\draw(g6)--(f1)--(g1)--(g3)--(g4)--(g5)--(g2)--(f2)--(g7);
\draw(g1)--(g6)--(g3);
\draw(g2)--(g7)--(g5);
\draw(t3)--(p1)--(p3)--(c1)--(c3)--(c5)--(c7)--(c6)--(c4)--(c2)--(c1);
\draw(p2)--(p3)--(p4);
\end{tikzpicture}
\caption{The graph $G_i$ corresponding to a variable $x_i$.}\label{graph}
\end{figure}
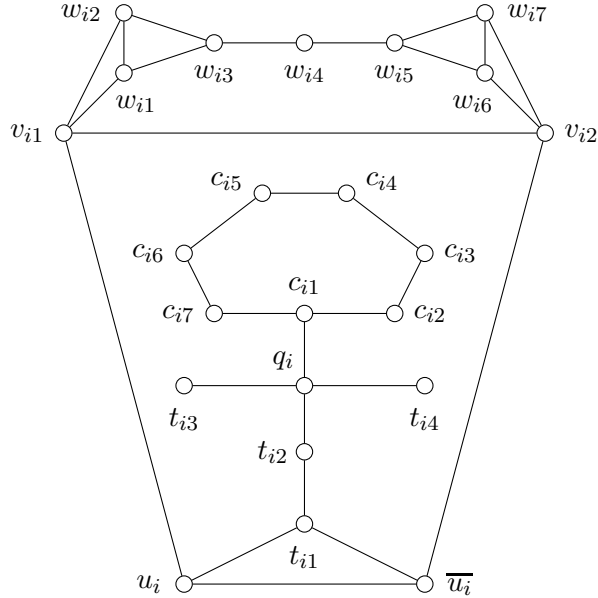

\proof
It is clear that the {\sc EOCD Problem} is in NP, since verifying that a given set of
vertices of a graph $G$ is an EOD set or an ECD set
can be done in polynomial time. We consider now a Boolean formula $\mathcal{F}$ {\s with variables $X=\{x_1,\dots,x_n\}$} {\iztok and} {\s clauses $C=\{c_1,\dots,c_m\}$}. Each clause contains three literals, each of which we shall denote
by $x_i$ for a positive literal, or by $\overline{x_i}$ for a negative one. From
the formula $\mathcal{F}$, we construct a graph $G_F$ in the following way.
For any variable $x_i\in X$, add to $G_F$ {\s the graph $G_i$ from} Fig.~\ref{graph}.
For each clause $c_i\in C$, we add a vertex $y_i$. Now, if a variable $x_i$ occurs
as a positive literal in a clause $c_j$, then add the edge $y_ju_i$, otherwise
(if a variable $x_i$ occurs as a negative literal in a clause $c_j$) add the edge $y_j\overline{u_i}$. Clearly, $G_F$ can be constructed in polynomial time.

We claim that $G_F$ is an EOCD graph if and only if there is a satisfying truth assignment for the $n$ variables in the Boolean formula $\mathcal{F}$, such that each clause has exactly one true literal, that is, if and only if $\mathcal{F}$ has
a one-in-three satisfying truth assignment.

Assume first that $\mathcal{F}$ has a one-in-three satisfying truth assignment. We construct two
sets $D$ and $P$ in the following way. Add to $D$ the vertices $q_i,c_{i1},c_{i4},c_{i5}$,
and to $P$ the vertices $q_i,c_{i3},c_{i6}$ for every $i\in [n]$. Now, if the
variable $x_i$ is assigned the value true, then we add to $D$ the vertices
$u_i,v_{i1},w_{i4},w_{i5}$, and to $P$ the vertices $u_i,w_{i3},w_{i7}$. On the other
hand, if $x_i$ is assigned the value false, then we add to {\iztok $D$} the vertices
$\overline{u_i},v_{i2},w_{i3},w_{i4}$, and to $P$ the vertices $\overline{u_i},w_{i2},w_{i5}$.
It is easy to see that $D\cap V(G_i)$ is an EOD set and $P\cap V(G_i)$ an ECD set of $G_i$. Moreover, since the
truth assignment has exactly one literal with value true, each vertex $y_j$, with
$j\in [m]$, is adjacent to exactly one vertex of $D$ and exactly one vertex of
$P$ (clearly both vertices coincide). Thus, $D$ is an EOD set and
$P$ is an ECD set in $G_F$ and, as a consequence, $G_F$ is an EOCD graph.

Conversely, assume that $G_F$ is an EOCD graph. Let $D$ be an EOD set and $P$ an ECD set in $G_F$. We next collect several facts regarding the sets $D$ and $P$.

\begin{itemize}
\item The vertex $q_i$ ($i\in [n]$) lies in $D\cap P$. Indeed, this fact follows because $q_i$ is adjacent to leaves $t_{i3}$ and $t_{i4}$.
\item The vertices $c_{i1}, c_{i4}, c_{i5}$ ($i\in [n]$) belong to $D$.
The vertex $c_{i1}$ belongs to $D$ because otherwise the $7$-cycle on vertices $c_{ij}$ cannot be efficiently open dominated. We then consequently see that also $c_{i4},c_{i5}\in D$.
\item The vertices $t_{i2},t_{i3},t_{i4},c_{i1}$ ($i\in [n]$) do not lie in $P$, and the vertices  $c_{i3},c_{i6}$ ($i\in [n]$) lie in $P$.
These facts follow immediately from the first point.
\item Either ($u_i\notin D$ and $\overline{u_i}\in D$) or ($u_i\in D$ and $\overline{u_i}\notin D$). Similarly, either ($u_i\notin P$ and
$\overline{u_i}\in P$) or ($u_i\in P$ and $\overline{u_i}\notin P$).
Indeed, since $q_i\in D\cap P$ and $c_{i1}\in D$, for every $i\in [n]$, the vertices
$t_{i1},t_{i2}\notin D\cup P$. Thus, $t_{i1}$ must be dominated either by $u_i$ or by $\overline{u_i}$ in $D$ and in $P$.
\item If $u_i\in P$ and $\overline{u_i}\notin P$, then $v_{i1}\notin P$ and every vertex $y_j$
such that the variable $x_i$ belongs to the clause $c_j$ does not belong to $P$. Moreover, to
efficiently dominate the vertices $v_{i2},w_{i1},\dots,w_{i7}$ we clearly have that
$w_{i3}\in P$ and exactly one vertex {\s of} the pair $w_{i6},w_{i7}$ belongs to $P$.
\item Analogously, if $u_i\notin P$ and $\overline{u_i}\in P$, then we obtain that
$v_{i2}\notin P$ and every vertex $y_j$ such that the variable $x_i$ belongs to the clause
$c_j$ does not belong to $P$. Also, $w_{i5}\in P$ and exactly one vertex {\s of} the pair
$w_{i1},w_{i2}$ belongs to $P$.
\item If $u_i\in D$ and $\overline{u_i}\notin D$, then either $v_{i1}\in D$ or there exists
a vertex $y_j\in D$ such that the variable $x_i$ appears as positive in the clause $c_j$. If
the latter happens ($y_j\in D$), then $v_{i1}\notin D$. It is straightforward to observe that,
in such a case, any subset of vertices {\igy of the set $\{v_{i2},w_{i1},\dots,w_{i7}\}$ does not efficiently open dominate the same set of vertices $\{v_{i2},w_{i1},\dots,w_{i7}\}$}, which is a contradiction. Thus, $y_j\notin D$ and
therefore $v_{i1}\in D$. We also observe that $w_{i4},w_{i5}\in D$.
\item Analogously to the last item, if $u_i\notin D$ and $\overline{u_i}\in D$, then $v_{i2},w_{i3},w_{i4}\in D$.
\end{itemize}

As a consequence of the {\s above facts}, we have that either $u_i,v_{i1}\in D$ or
$\overline{u_i},v_{i2}\in D$, and either $u_i\in P$ or $\overline{u_i}\in P$. Now, we say that a subgraph $G_i$ of $G_F$, corresponding to a variable $x_i$, is \textit{nice} if either $u_i\in D\cap P$ or $\overline{u_i}\in D\cap P$. Assume that there exists $G_i$ which is not nice, \emph{i.e.},
$u_i,\overline{u_i}\notin D\cap P$. Hence, either $u_i\in D$ and $\overline{u_i}\in P$ or $u_i\in P$ and
$\overline{u_i}\in D$. Consider a clause $c_j$ such that $x_i\in c_j$. Hence, $y_j$ is dominated either by
$u_i$ or by $\overline{u_i}$ from $G_i$, which means either by $D$ or by $P$, say $D$. Therefore, there must exist another variable
$x_{\ell}\in c_j$, such that $y_j$ is dominated also by $P$ and not by $D$. This implies that $G_{\ell}$ is not nice as well.
Notice that for the third literal $x_k\in c_j$, $G_k$ must be nice. In general, for every clause $c_j$, either all three corresponding graphs are nice,
or exactly one is nice and two are not. Moreover, if the later is true, then $y_j$ is not dominated from $D$ and from $P$ by the nice subgraph.

Let $P_i=P\cap V(G_i)$ ($i\in [n]$). For every not nice graph $G_i$ we exchange some vertices of $P_i$ as follows. If $u_i\in P_i$, then
$P'_i=(P_i-\{u_i,w_{i3},w_{i6},w_{i7}\})\cup \{{\overline{u_i}},w_{i2},w_{i5}\}$, and if $\overline{u_i}\in P$, then
$P'_i=(P_i-\{{\overline{u_i}},w_{i5},w_{i1},w_{i2}\})\cup \{u_i,w_{i3},w_{i7}\}$. If $G_i$ is nice, then $P'_i=P_i$.
We claim that $P'=\cup_{i=1}^n P'_i$ is an ECD set, such that together with $D$ every subgraph $G_i$ is nice. Clearly $P'_i$ is an ECD set for $G_i$ by the items above. If some $y_j$ was dominated by a vertex $u_i$ (or by $\overline{u_i}$) which was in $P_i$ but not now in $P'_i$, then $y_j$ is now dominated either by
$u_{\ell}$ or by $\overline{u_{\ell}}$, where $x_i$ and $x_{\ell}$ are those variables from the clause $c_j$, for which $G_i$ and $G_{\ell}$ were not nice. Thus, $P'$ is an ECD set. Moreover, the EOD set $D$ and the ECD set $P'$ lead to {\s the fact} that every $G_i$ is nice.
Since $D$ is an EOD set and $P'$ is an ECD set, then every vertex $y_j$ corresponding
to a clause is adjacent to exactly one vertex $u_i$ or $\overline{u_i}$ of $G_i$.
Now, if $u_i\in D\cap P'$, then we set the variable $x_i$ as true, {\s otherwise (if
$\overline{u_i}\in D\cap P'$) set $x_i$ as false}. It clearly follows
that such an assignment is a truth assignment in exactly one literal in every clause for $\mathcal{F}$ and the proof is complete.
\qed

In view of Theorem~\ref{thm:NP-complete}, it is reasonable to try to find some special classes of graphs for which the {\sc EOCD Problem} is polynomial. Simple examples are provided by the paths $P_{n}$ which are EOCD graphs if and only if $n\not\equiv 1\pmod{4}$ and the cycles $C_{n}$ which are EOCD graphs if and only if $n\equiv 0\pmod{12}$. Note also that a complete bipartite graph $K_{r,t}$ is an EOCD graph if and only if $r=1$ or $t=1$. {\igy Moreover, the hypercube $Q_n$ is an EOCD graph if and only if $n=1$. Indeed, suppose that $Q_n$, $n\ge 1$, is an EOCD graph. As $Q_n$ is $n$-regular, its order must be divisible by $n$ (because it admits an
EOD set) as well as by $n+1$ (since it admits an ECD set). Since the order of $Q_n$ is $2^n$, this is only possible if $n=1$.} 

We end this section with a discussion on extreme cases with respect to the structure of EOCD graphs as described in Section~\ref{sec:structure}.

\begin{theorem}
If $G$ is a graph on $n$ vertices and $m$ edges, then it can be decided in $O(nm)$ time whether $G$ is an EOCD graph with empty $P-D$.
\end{theorem}

\proof
Let $G$ be a graph. Clearly, components which are isomorphic to $K_2$ (if they exists) do not influence the fact that $G$ is an EOCD graph or not.
Hence we may restrict our attention to the case when $G$ has no components isomorphic to $K_2$. If there exists no degree {\s $1$} vertex,
then by Proposition~\ref{prp:P-D:empty}, $G$ is not an EOCD graph with empty $P-D$. Let $P$ be the set of all support vertices of
degree one vertices. For every support from $P$ choose exactly one neighbor of degree {\s $1$} and let $D$ be a set containing $P$ as well as the
chosen vertices of degree {\s $1$}. By Proposition~\ref{prp:P-D:empty} one only needs to check if $D$ and $P$ are an
EOD set and an ECD set of $G$, respectively. Even more, it is clear that $P$ is an ECD set in $G$ if and only if $D$ is an EOD set
of $G$. Hence it is enough to check
whether the union of closed neighborhoods centered in $P$ covers $V(G)$ and whether these closed neighborhoods have pairwise empty intersection. The first task can be clearly done in $O(m)$ time.
For the second task it suffices to check if the distance between any two different vertices from $P$ is at least 3. Clearly, this can be done in {\s time} $O(mn)$,
if we start the BFS algorithm in {\s an arbitrary} vertex of $P$.
\qed

We end the section with a question about the other extremal case.

\begin{problem}
{\s Can} it be checked in polynomial time whether $G$ is an EOCD graph with empty $D\cap P$?
\end{problem}


\section{EOCD trees}
\label{sec:trees}

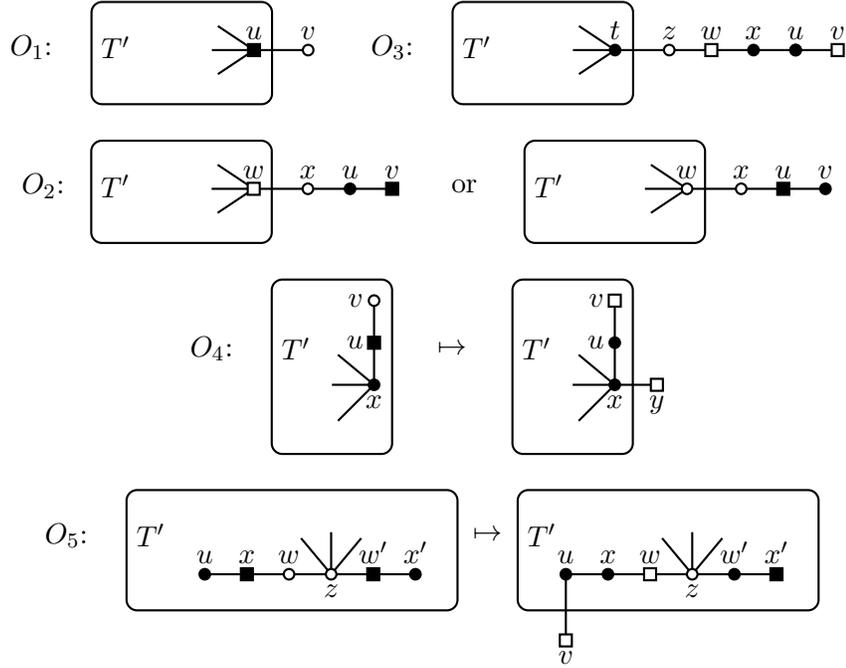
\begin{figure}[ht!]
\begin{center}
\begin{tikzpicture}[scale=.8,style=thick,x=1cm,y=1cm]
\def\vr{2.5pt} 
\path (3.7,0.7) coordinate (v);
\path (4.6,0.7) coordinate (v1);
\path (3.1,1.1) coordinate (a);
\path (3.0,0.7) coordinate (b);
\path (3.1,0.3) coordinate (c);
%
\draw (v) -- (v1);
\draw (v) -- (a);
\draw (v) -- (b);
\draw (v) -- (c);

\draw (v1) [fill=white] circle (\vr);
\draw [fill=black] (3.6,0.8) rectangle (3.8,0.6);
\draw [rounded corners] (1,-.2) rectangle (4,1.5);
\draw (1.4,0.75) node {$T'$};
\draw (0,.75) node {$O_{1}$:};
\draw[anchor = south] (v) node {$u$};
\draw[anchor = south] (v1) node {$v$};

\path (9.7,0.7) coordinate (v);
\path (10.6,0.7) coordinate (v1);
\path (11.3,0.7) coordinate (v2);
\path (12.0,0.7) coordinate (v3);
\path (12.7,0.7) coordinate (v4);
\path (13.4,0.7) coordinate (w);
\path (9.1,1.1) coordinate (a);
\path (9.0,0.7) coordinate (b);
\path (9.1,0.3) coordinate (c);
%
\draw (v) -- (v1);
\draw (v1) -- (v2);
\draw (v2) -- (v3);
\draw (v3) -- (v4);
\draw (v4) -- (w);
\draw (v) -- (a);
\draw (v) -- (b);
\draw (v) -- (c);

\draw (v) [fill=black] circle (\vr);
\draw (v1) [fill=white] circle (\vr);
\draw (v4) [fill=black] circle (\vr);
\draw (v3) [fill=black] circle (\vr);
\draw (v4) [fill=black] circle (\vr);
\draw (w) [fill=white] circle (\vr);
\draw [fill=white] (11.2,0.8) rectangle (11.4,0.6);
\draw [fill=white] (13.3,0.8) rectangle (13.5,0.6);
\draw [rounded corners] (7,-.2) rectangle (10,1.5);
\draw (7.4,0.75) node {$T'$};
\draw (6,.75) node {$O_{3}$:};
\draw[anchor = south] (v) node {$t$};
\draw[anchor = south] (v1) node {$z$};
\draw[anchor = south] (v2) node {$w$};
\draw[anchor = south] (v3) node {$x$};
\draw[anchor = south] (v4) node {$u$};
\draw[anchor = south] (w) node {$v$};

\end{tikzpicture}
\end{center}

\begin{center}
\begin{tikzpicture}[scale=.8,style=thick,x=1cm,y=1cm]
\def\vr{2.5pt} 
\path (4.5,0.7) coordinate (u);
\path (5.4,0.7) coordinate (u1);
\path (6.1,0.7) coordinate (u2);
\path (6.8,0.7) coordinate (u3);
\path (3.9,1.1) coordinate (a2);
\path (3.8,0.7) coordinate (b2);
\path (3.9,0.3) coordinate (c2);
%
\draw (u) -- (u1);
\draw (u1) -- (u2);
\draw (u2) -- (u3);
\draw (u) -- (a2);
\draw (u) -- (b2);
\draw (u) -- (c2);

\draw (u2) [fill=black] circle (\vr);
\draw (u1) [fill=white] circle (\vr);
\draw (u3) [fill=black] circle (\vr);
\draw [fill=black] (6.7,0.8) rectangle (6.9,0.6);
\draw [fill=white] (4.4,0.8) rectangle (4.6,0.6);
\draw [rounded corners] (1.8,-.2) rectangle (4.8,1.5);
\draw (2.2,0.75) node {$T'$};
\draw (1,.75) node {$O_{2}$:};
\draw[anchor = south] (u) node {$w$};
\draw[anchor = south] (u1) node {$x$};
\draw[anchor = south] (u2) node {$u$};
\draw[anchor = south] (u3) node {$v$};

\draw (8,.75) node {or};

\path (11.7,0.7) coordinate (w);
\path (12.6,0.7) coordinate (w1);
\path (13.3,0.7) coordinate (w2);
\path (14.0,0.7) coordinate (w3);
\path (11.1,1.1) coordinate (a1);
\path (11.0,0.7) coordinate (b1);
\path (11.1,0.3) coordinate (c1);
%
\draw (w) -- (w1);
\draw (w1) -- (w2);
\draw (w2) -- (w3);
\draw (w) -- (a1);
\draw (w) -- (b1);
\draw (w) -- (c1);

\draw (w) [fill=white] circle (\vr);
\draw (w1) [fill=white] circle (\vr);
\draw (w3) [fill=black] circle (\vr);
\draw [fill=black] (13.2,0.8) rectangle (13.4,0.6);
\draw [rounded corners] (9,-.2) rectangle (12,1.5);
\draw (9.4,0.75) node {$T'$};
\draw[anchor = south] (w) node {$w$};
\draw[anchor = south] (w1) node {$x$};
\draw[anchor = south] (w2) node {$u$};
\draw[anchor = south] (w3) node {$v$};

\end{tikzpicture}
\end{center}

\begin{center}
\begin{tikzpicture}[scale=.8,style=thick,x=1cm,y=1cm]
\def\vr{2.5pt} 
\path (1.7,1.15) coordinate (v);
\path (1.7,1.85) coordinate (v1);
\path (1.7,2.55) coordinate (v2);
\path (1.1,.55) coordinate (a);
\path (1.0,1.15) coordinate (b);
\path (1.1,1.65) coordinate (c);
%
\draw (v) -- (v1);
\draw (v1) -- (v2);
\draw (v) -- (a);
\draw (v) -- (b);
\draw (v) -- (c);

\draw (v) [fill=black] circle (\vr);
\draw (v2) [fill=white] circle (\vr);
\draw [fill=black] (1.6,1.95) rectangle (1.8,1.75);
\draw [rounded corners] (0,0) rectangle (2,2.9);
\draw (.4,1.75) node {$T'$};
\draw (-1,1.75) node {$O_4$:};
\draw[anchor = north] (v) node {$x$};
\draw[anchor = east] (v1) node {$u$};
\draw[anchor = east] (v2) node {$v$};

\draw (3,1.75) node {$\mapsto$};
\path (5.7,1.15) coordinate (w);
\path (5.7,1.85) coordinate (w1);
\path (5.7,2.55) coordinate (w2);
\path (5.1,.55) coordinate (a1);
\path (5.0,1.15) coordinate (b1);
\path (5.1,1.65) coordinate (c1);
\path (6.4,1.15) coordinate (y);
%
\draw (w) -- (w1);
\draw (w1) -- (w2);
\draw (w) -- (a1);
\draw (w) -- (b1);
\draw (w) -- (c1);
\draw (w) -- (y);

\draw (w) [fill=black] circle (\vr);
\draw (w1) [fill=black] circle (\vr);
\draw [fill=white] (5.6,2.65) rectangle (5.8,2.45);
\draw [fill=white] (6.3,1.05) rectangle (6.5,1.25);
\draw [rounded corners] (4,0) rectangle (6,2.9);
\draw (4.4,1.75) node {$T'$};
\draw[anchor = north] (w) node {$x$};
\draw[anchor = east] (w1) node {$u$};
\draw[anchor = east] (w2) node {$v$};
\draw[anchor = north] (y) node {$y$};
\end{tikzpicture}
\end{center}

\begin{center}
\begin{tikzpicture}[scale=.8,style=thick,x=1cm,y=1cm]
\def\vr{2.5pt} 
\path (1.3,1.1) coordinate (u);
\path (2.0,1.1) coordinate (x);
\path (2.7,1.1) coordinate (w);
\path (2.9,1.7) coordinate (a);
\path (3.4,1.8) coordinate (b);
\path (3.9,1.7) coordinate (c);
\path (3.4,1.1) coordinate (z);
\path (4.1,1.1) coordinate (t);
\path (4.8,1.1) coordinate (r);

%
\draw (u) -- (x) -- (w) -- (z) -- (t) -- (r);
\draw (z) -- (a);
\draw (z) -- (b);
\draw (z) -- (c);

\draw (u) [fill=black] circle (\vr);
\draw (r) [fill=black] circle (\vr);
\draw (w) [fill=white] circle (\vr);
\draw (z) [fill=white] circle (\vr);
\draw [fill=black] (1.9,1.2) rectangle (2.1,1.0);
\draw [fill=black] (4.0,1.2) rectangle (4.2,1.0);
\draw [rounded corners] (0,0.5) rectangle (5.5,2.5);
\draw (.4,1.75) node {$T'$};
\draw (-1,1.75) node {$O_5$:};
\draw[anchor = south] (u) node {$u$};
\draw[anchor = south] (x) node {$x$};
\draw[anchor = south] (w) node {$w$};
\draw[anchor = north] (z) node {$z$};
\draw[anchor = south] (t) node {$w'$};
\draw[anchor = south] (r) node {$x'$};

\draw (6,1.75) node {$\mapsto$};
\path (7.3,1.1) coordinate (u1);
\path (8.0,1.1) coordinate (x1);
\path (8.7,1.1) coordinate (w1);
\path (8.9,1.7) coordinate (a1);
\path (9.4,1.8) coordinate (b1);
\path (9.9,1.7) coordinate (c1);
\path (9.4,1.1) coordinate (z1);
\path (10.1,1.1) coordinate (t1);
\path (10.8,1.1) coordinate (r1);
\path (7.3,0) coordinate (v);

%
\draw (u1) -- (x1) -- (w1) -- (z1) -- (t1) -- (r1);
\draw (z1) -- (a1);
\draw (z1) -- (b1);
\draw (z1) -- (c1);
\draw (u1) -- (v);

\draw (u1) [fill=black] circle (\vr);
\draw (x1) [fill=black] circle (\vr);
\draw (t1) [fill=black] circle (\vr);
\draw (z1) [fill=white] circle (\vr);
\draw [fill=white] (8.6,1.2) rectangle (8.8,1.0);
\draw [fill=white] (7.2,-0.1) rectangle (7.4,0.1);
\draw [fill=black] (10.7,1.2) rectangle (10.9,1.0);
\draw [rounded corners] (6.5,0.5) rectangle (11.5,2.5);
\draw (6.9,1.75) node {$T'$};
\draw[anchor = south] (u1) node {$u$};
\draw[anchor = south] (x1) node {$x$};
\draw[anchor = south] (w1) node {$w$};
\draw[anchor = north] (z1) node {$z$};
\draw[anchor = south] (t1) node {$w'$};
\draw[anchor = south] (r1) node {$x'$};
\draw[anchor = north] (v) node {$v$};
\end{tikzpicture}
\end{center}

\caption{Operations $O_{1}-O_{5}$.}
\label{operation123}
\end{figure}

Let $T'$ be an EOCD tree with an EOD set $D'$ and an ECD set $P'$. We now introduce five operations that construct larger EOCD trees from $T'$. In the main theorem of this section we will then prove that these operations are characteristic for EOCD trees. The operations are illustrated in Fig.~\ref{operation123} where we use the convention introduced in Section~\ref{sec:structure}: a vertex from $D\cap P$ is black squared, a vertex from $P-D$ is white squared, a vertex from $D-P$ is black circled, and the remaining vertices are white circled.

\begin{itemize}
\item[$(O_1)$] For $u\in D'\cap P'$ we obtain $T$ from $T'$ by adding a vertex $v$ and edge $uv$.
\item[$(O_2)$] For $w\notin D'$ we obtain $T$ from $T'$ by
adding a path $xuv$ and edge $wx$.
\item[$(O_3)$] For $t\in D'-P'$ we obtain $T$ from $T'$ by
adding a path $zwxuv$ and edge $tz$.
\item[$(O_4)$] For a path $vux$ with $\deg(v)=1$, $\deg(u)=2$, $u,x\in
D'$, and $u\in P'$, we obtain $T$ from $T'$ by
adding a vertex $y$ and edge $yx$.
\item[$(O_5)$] For a path $uxwzw'x'$ with $\deg(u)=\deg(x')=1$, $\deg(x)=\deg(w)=\deg(w')=2$, $u,x,w',x'\in D'$, and $x,w'\in P'$, we obtain $T$ from $T'$
by adding a vertex $v$ and edge $uv$.
\end{itemize}

The main difference between the{\s se} five operations is that for $O_1$ the original EOD set and ECD set do not change, for $O_2$ and $O_3$ we add some vertices to the EOD set and to the ECD set, while for $O_4$ and $O_5$ the EOD set remains the same and we exchange some vertices in the ECD set.

\begin{theorem}\label{tree}
A tree $T$ is an EOCD graph if and only if $T$ can be obtained
from $K_{2}$ by a sequence of operations $O_1-O_5$.
\end{theorem}

\proof
Assume first that $T$ is a tree obtained from $K_{2}$ by a sequence
of operations $O_1-O_5$. We will show that $T$ is an EOCD tree by induction
on the length $k$ of the mentioned sequence. If $k=0$, then $T\cong K_{2}$
which is an EOCD graph. Let now $k>0$ and let $T'$ be a tree
obtained from $K_{2}$ by {\igy using} the same sequence {\igy as for $T$, but without including the} {\iztok last} step. By the
induction hypothesis, $T'$ is an EOCD tree with an EOD set $D'$ and an ECD set $P'$. If $T$ is obtained from $T'$
by operation $O_1$, then clearly $T$ is an EOCD tree for $D=D'$ and $P=P'$ (see the upper left diagram of Fig.~\ref{operation123}).
If $T$ is obtained from $T'$ by operation $O_2$, then $T$ is an EOCD tree where $D=D'\cup \{u,v\}$. The set $P$ depends
whether $w$ is in $P'$ or not. If $w\in P'$, then $P=P'\cup \{v\}$, and if $w\notin P'$, then $P=P'\cup \{u\}$ (see the diagrams of the second line of Fig.~\ref{operation123}).
Suppose now that we apply operation $O_3$ on $T'$ to get
$T$. Again it is straightforward to see that $T$ is an EOCD graph for $D=D'\cup \{u,x\}$ and $P=P'\cup \{v,w\}$ (see the upper right diagram of Fig.~\ref{operation123}). If operation $O_4$
is applied to get $T$ from $T'$, then we set $D=D'$ and $P=(P'-\{u\})\cup \{v,y\}$ and $T$ is an EOCD tree again (see the diagram in the third line of Fig.~\ref{operation123}). Finally,
if $T$ is obtained from $T'$ by operation $O_5$, then it is not hard to
see that $T$ is an EOCD tree for $D=D'$ and $P=(P'-\{x,w'\})\cup \{v,x',w\}$ (see the lower diagram of Fig.~\ref{operation123}).

To prove the converse, let $T$ be an EOCD tree with an EOD set $D$ and an ECD set $P$. Let $r\in V(T)$ be a vertex of $T$ and consider $T$ as a rooted tree with the root $r$. Let $v$ be a vertex of degree {\s $1$} of $T$ that is at the maximum distance from $r$ and let $u$ be the support vertex of $v$. Clearly $u\in D$, while either {\igy $u\in P$ or $v\in P$}. We  call a neighbor $y$ of $x$ a \textit{down}- (resp.\ \textit{up}-) \textit{neighbor} of $x$ if $y$ is further (resp.\ closer) {\s from} $r$ than $x$.  We proceed by induction on the number of vertices of $T$. Clearly, $K_2$ is the smallest EOCD tree, hence the base of the induction. We distinguish the following cases.

\medskip\noindent
{\bf Case 1:} $v\notin P$ and $v\notin D$. \\
In this case $u\in P\cap D$. We obtain a tree $T'$ from $T$ by deleting $v$. Clearly $T'$
is an EOCD tree with $D'=D$ and $P'=P$. By the induction
hypothesis $T'$ can be buil{\igy t} from $K_{2}$ by a sequence of operations $O_1-O_5$. If we add the operation $O_1$ at the end of this sequence, then we obtain $T$ from $K_{2}$ by a sequence of operations $O_1-O_5$.

\medskip\noindent
{\bf Case 2:} $v\notin P$ and $v\in D$. \\
Then $u\in P\cap D$. If $\deg(u)=1$, then $T\cong K_{2}$ and we are done. So, let $\deg(u)>1$.
If $u$ is the support for more degree {\s $1$} vertices than $v$, then we have
Case 1. (Notice that the same does occur when $u=r$.) Thus let $\deg(u)=2$.
Let $x$ be the up-neighbor of $u$. If $\deg(x)>2$, then $x$ has a down-neighbor
$y$ different from $u$. If $\deg(y)=1$, then we have a contradiction with
$P$ being an ECD set of $T$, since $u\in P$ implies that $y$ and $x$
cannot be in $P$ and therefore $y$ is {\s neither} dominated by $P$ nor $y\in P$. So $\deg(y)>1$
and let $y'$ be a down-neighbor of $y$. Clearly, $\deg(y')=1$ by the choice of $v$. This yields a contradiction with $D$ being an
EOD set of $T$, since $y$ cannot be in $D$ because $x$ is already dominated
by $u\in D$. Thus, $\deg(x)=2$ and let $w$ be the up-neighbor of $x$ or
the other down-neighbor when $x=r$. By the choice of $v$, $x$ must be
different from $r$ or we obtain the same problems as for $\deg(x)>2$.
Since $x$ is the neighbor of $u\in D\cap P$, we have that $w\notin D$ and $w\notin P$. Let $T'$ be the tree obtained from $T$ by deleting
vertices $v, u, x$. Then $T'$ is an EOCD tree with $D'=D-\{u,v\}$ and $P'=P-\{u\}$. By the induction hypothesis,
$T'$ can be buil{\igy t} from $K_{2}$ by a sequence of operations $O_1-O_5$. Adding operation $O_2$ {\igy at} the end of this sequence we obtain $T$ from $K_{2}$ by a sequence of operations $O_1-O_5$.

\medskip\noindent
{\bf Case 3:} $v\in P\cap D$. \\
{\iztok If $\deg(u)=1$, then $T\cong K_{2}$ and we are done.} If $\deg(u)>2$, then we have a
contradiction with $P$ being an ECD set. Thus $\deg(u)=2$. Also
notice that $u\neq r$, since otherwise $T$ would be isomorphic to {\iztok $P_3$},
which is not possible with $v\in P$. Let $x$ be the up-neighbor of $u$.
Clearly $x\notin D\cup P$. If $x$ would have a down-neighbor different from
$u$ or if $x=r$, then we have a contradiction with $D$ being an EOD set
of $T$. Thus, also $\deg(x)=2$ and let $w$ be the up-neighbor of $x$.
Notice that now $w$ must be in $P$ to dominate $x$, but again $w\notin D$.
Let $T'$ be the tree obtained from $T$ by deleting vertices $v, u, x$. Clearly $T'$ is an EOCD tree with $D'=D-\{u,v\}$ and $P'=P-\{v\}$. By the induction hypothesis $T'$ can be buil{\igy t}
from $K_{2}$ by a sequence of operations $O_1-O_5$. Adding the corresponding operation $O_2$ {\igy at} the end of this sequence we obtain $T$ from $K_{2}$ {\s as desired}.

\medskip\noindent
{\bf Case 4:} $v\in P$ and $v\notin D$. \\
In this case $u\notin P$.
If $u=r$, then we have a contradiction with $P$ being an ECD set when $\deg(u)>1$
and with $D$ being an OED set if $\deg(u){\iztok =}1$. So we may assume that $u\neq r$.
Clearly $\deg(u)=2$, otherwise we have a
contradiction again with $P$ being an ECD set of $T$ {\iztok and by the choice of $v$}. Let $x$ be the
up-neighbor of $u$. Since $v\notin D$ and $v\in P$, we have that $x\in D$
and $x\notin P$, respectively. The only second down-neighbor of $x$ is $v$,
otherwise we have a contradiction with $D$ being an EOD set for $T$ {\igy according to} the
choice of $v$. Suppose that $x$ has a down-neighbor $y$ of degree {\s $1$}.
Clearly, $v\in P$ implies $x\notin P$ and therefore $y\in P$ and $y$ is {\igy the} unique
down-neighbor of $x$ of degree {\s $1$}. Thus $vuxy$ is a path. Let $T'$
be a tree obtained from $T$ by deleting {\igy the} vertex $y$. Clearly $T'$ is
an EOCD tree with $D'=D$ and $P'=\left( P-\{v,y\}\right)
\cup \{u\}$. By the induction hypothesis $T'$ can be buil{\igy t} from $K_{2}$ by a sequence of operations $O_1-O_5$. If we add the operation $O_4$ at the end of this
sequence, then we obtain $T$ from $K_{2}$ by a sequence of operations
$O_1-O_5$.

Suppose now that $x$ has no down-neighbor of degree {\s $1$}. If $x=r$, then we
have a contradiction with $P$ being an ECD set for $T$. Hence $x\neq r$
and $\deg(x)=2$ holds. Let $w$ be the up-neighbor of $x$. Clearly, $w\in
P $ and $w\notin D$. If $\deg(w)\geq 3$, then $w$ has a down-neighbor $x'$ other than $x$. To dominate $x'$ from $D$, the vertex $x'$ must have a down-neighbor $u'$ which is in $D$ and the
same holds for $u'$, which must have a down-neighbor $v'$
which is also in $D$. Moreover, to dominate $x',u'$ and $v'$ from $P$ exactly once, also $v'\in P$ {\igy holds}. Notice that $\deg(x')=2$ {\igy according to that $D$ is an EOD set of $T$, and $\deg(u')=2$ since $P$ is an ECD set} of $T$. The situation for $v'$ is now as in Case 3 and we are done if $\deg(w)\geq 3$.

Thus, from now on, we consider $\deg(w)=2$ and let $z$ be the up-neighbor of $w$
(or down-neighbor if $w=r$). Again $z\notin P$ since $w\in P$, and $z\notin D$ since
$x\in D$.  We consider the following subcases.

\medskip\noindent
{\bf Subcase 4.1:} $\deg(z)\geq 3$. \\
{\s Let} $w'\ne w$ be a down-neighbor
of $z$. Since $z$ is dominated from $P$ by $w$, $w'$ is not in $P$
and therefore, $w'$ must have a down-neighbor $x'$ which is
in $P$. Also, $x'$ cannot have a second down-neighbor by the choice of $v$ and the structure of $P$. We consider two possibilities regarding the vertex $w'$.

\medskip\noindent
{\bf Subcase 4.1.1:} $w'\notin D$. \\
{\s In this subcase} $w'$ must be dominated by a
down-neighbor in $D$. If $x'\notin D$, then we have a contradiction since $x'$ has no second-down
neighbor and $D$ is an EOD set. Hence, $x'\in D$ and $x'$ must have a down-neighbor $u'\in D$ for $x'$ to be
dominated by $D$. Clearly $u'\notin P$. If $x'$ has another down-neighbor $u''$, then we obtain a
tree $T'$ from $T$ by deleting $v$. Clearly $T'$ is an EOCD tree with $D'=D$ and $P'=P$. By the induction
hypothesis, $T'$ can be built from $K_{2}$ by a sequence of operations $O_1-O_5$. If we add the operation $O_1$ at the end
of this sequence, then we obtain $T$ from $K_{2}$ by a sequence of operations $O_1-O_5$.
So we may assume that $\deg(x')=2$. If $\deg(w')=2$, then let $T'$ be a tree obtained from $T$ by deleting vertices
$u', x', w'$. Clearly $T'$ is an EOCD tree with $D'=D-\{u',x'\}$ and $P'=P-\{x'\}$. By the induction hypothesis, $T'$ can be built
from $K_{2}$ by a sequence of operations $O_1-O_5$. Adding the corresponding operation $O_2$ at the end of this sequence
we obtain $T$ from $K_{2}$ as desired. On the other hand, if $w'$ has a down-neighbor $x''$ other than $x'$,
then $x''$ is not in $D$ and not in $P$, since $w'$ is already dominated by $x'$ in both $P$ and $D$.
Moreover, $x''$ must be dominated by its down-neighbor $u''$ in both $P$ and $D$. Furthermore, $u''$ is dominated by its down-neighbor $v''$ in $D$.
If $\deg(u'')>2$, then we have Case 1. So let $\deg(u'')=2$. If $\deg(x'')>2$, we have a contradiction with $D$ being an EOD set of
$T$ or by the choice of $v$. Hence, $\deg(x'')=2$ and {\igy we proceed like in} Case 2 for $v'',u'', x''$.

\medskip\noindent
{\bf Subcase 4.1.2:}  $w'\in D$. \\
{\s Since $z\notin D$, also in this subcase $w'$ must have a down-neighbor in $D$}. Suppose first that $x'\in D$ (and recall that $x'\in P$). If $x'$ has a down-neighbor $u'$, then $\deg(u')=1$ by the choice of $v$ and since $P$ is an ECD set. Let $T'$ be a tree obtained from $T$
by deleting $u'$. Clearly $T'$ is an EOCD tree with $D'=D$ and $P'=P$. By the induction hypothesis $T'$ can be built from $K_{2}$ by a
sequence of operations $O_1-O_5$, and attaching operation $O_1$ to this sequence we obtain $T$ from $K_{2}$ as desired. Thus, we may assume that $\deg(x')=1$. Observe that $\deg(w')=2$, otherwise we have a contradiction with the choice of $v$, since $P$ is an ECD set
and since $D$ is an EOD set. Let $T'$ be a tree obtained from $T$ by deleting $v$. Clearly $T'$ is an EOCD tree with
$D'=D$ and $P'=(P-\{x',w{\iztok ,v}\})\cup \{x,w'\}$. By the induction hypothesis, $T'$ can be built from $K_{2}$ by a
sequence of operations $O_1-O_5$. Now, adding operation $O_5$ at the end of such sequence produces our desired result.
Next, let $x'\notin D$. Clearly $\delta (x')=1$, since $D$ is an EOD set and by the choice of $v$. Let $w$ now be dominated by $x''\in D$.
To dominate $x''$ from $P$, let $u''$ be its down-neighbor. Also, $\delta (u'')=1$ since $D$ is an EOD set and by the choice of $v$.
Observe that $\delta (x'')=2$ since any other down-neighbor $u'''$ of $x''$ would required a down-neighbor $v'$ in $P$, which is not possible
since $D$ is an EOD set and by the choice of $v$. Let $T'$
be a tree obtained from $T$ by deleting vertex $x'$. Clearly $T'$ is
an EOCD tree with $D'=D$ and $P'=\left( P-\{x',u''\}\right)\cup \{x''\}$. By the induction hypothesis, $T'$ can be built from
$K_{2}$ by a sequence of operations $O_1-O_5$. If we add operation $O_4$ at the end of this
sequence, then we obtain $T$ from $K_{2}$ by a sequence of operations $O_1-O_5$.

\medskip\noindent
{\bf Subcase 4.2:} $\deg(z)=2$. \\
Let $T'$ be a tree obtained from $T$ by
deleting $v, u, x, w, z$. Clearly $T'$ is an EOCD tree
with $D'=D-\{u,x\}$ and $P'=P-\{v,w\}$. Applying the induction hypothesis once more and ending with an additional operation $O_3$, we {\s again} obtain $T$ from $K_{2}$ {\s as desired and we are done}.
\qed

It is not obvious that all the five operations are necessary to characterize EOCD trees. To see that this is the case, note first that $P_{3}$ can be obtained from $K_{2}$ only by operation $O_1$ and that the sequence of operations $O_1,O_4$ is unique for $P_{4}$. Similarly, the sequence of operations $O_1,O_2$ is unique for $P_{6}$. {\s To infer that operations $O_3$ and $O_5$ are also indispensable, consider the following more elaborate examples}.

Let $T$ be the tree obtained from $K_{1,3}$ by subdividing {\iztok one of its edges with five vertices and each of the other two edges with eight vertices. A short analysis reveals that $T$ is an EOCD tree where the vertex of degree $3$ must be in $D\cap P$ and that its neighbor on the shortest leg must be in $D$. After this observation,}
operation $O_3$ cannot be avoided when constructing $T$ in view of Theorem~\ref{tree}. {\s For} operation $O_5$, let $P_{22}^{+}$ be the graph obtained from the path on $22$ vertices $v_{1},\ldots ,v_{22}$, by adding vertices $u,w,x,y$ and edges $v_{5}u,uw,v_{18}x,xy$. One can observe that $P_{22}^{+}$ is an EOCD tree with a unique EOD set $D$ and a unique ECD set $P$. {\s From} here it can be concluded that operation $O_5$ is needed to build $P_{22}^{+}$ from $K_2$ in view of Theorem~\ref{tree}. We leave the details to the reader.

\section{EOCD Sierpi\'nski graphs}

The Sierpi\'nski graphs $S_p^n$ were introduced in~\cite{klavzar-1997} and afterwards investigated from many different aspects. Here we only mention recent studies of Sierpi\'nski graphs related to codes and domination~\cite{dorbec-2014, gravier-2013, lin-2013}, their shortest paths~\cite{hinz-2014, xue-2014}, and an appealing generalization of Sierpi\'nski graphs due to Hasunuma~\cite{hasunuma-2015} that in turn extends several known results about Sierpi\'nski graphs. For the additional vast bibliography on these graphs we refer to~\cite{hasunuma-2015}.

The {\em Sierpi\'nski graphs} $S_p^n$, $p\ge 1$, $n\ge 0$, are defined as follows. $S_p^0 = K_1$ for any $p$. For $n\ge 1$, the vertex set of $S_p^n$ is $[p]_0^n$, we shall denote its elements by $s=s_n\ldots s_1$. Vertices $s_n\ldots s_1$ and $t_n\ldots t_1$ are adjacent if and only if there exists a $\delta\in [n]$ such that
\begin{itemize}
\item[(i)] $s_d = t_d$, for $d\in [n]-[\delta]$;
\item[(ii)] $s_\delta \ne t_\delta$; and
\item[(iii)] $s_d=t_\delta$ and $t_d=s_\delta$ for $d\in [\delta-1]$.
\end{itemize}
Note that $S_1^n \cong K_1$ ($n\ge 1$), $S_2^n \cong P_{2^n}$ ($n\ge 1$), and $S_p^1 \cong K_p$ ($p\ge 1$). Hence, for our purposes we may restrict the attention to the Sierpi\'nski graphs $S_p^n$ with $p\ge 3$ and $n\ge 2$.

The edge set of $S_p^n$ can be equivalently defined recursively as
$$E(S_p^{n}) = \{ \{is,it\}:\ i\in [p]_0\,, \{s,t\}\in E(S_p^{n-1} )\} \cup
\{ \{ij^{n-1}, ji^{n-1}\} \ |\ i,j\in [p]_0\,, i\ne j \}\,.$$
This implies that $S_p^{n}$ can be constructed from $p$ copies of $S_p^{n-1}$ as follows. For each $j\in [p]_0$ concatenate $j$ to the left of the vertices in a copy of $S_p^{n-1}$ and denote the obtained graph with $jS_p^{n-1}$. Then for each $i\neq j$ join copies $iS_p^{n-1}$ and $jS_p^{n-1}$ by the single edge $e_{ij}^{(n)} = \{ij^{n-1}, ji^{n-1}\}$.

If $1\le d < n$ and $\underline{s}\in [p]_0^{d}$, then the subgraph of $S_p^n$ induced by the vertices whose labels begin with $\underline{s}$ is isomorphic to $S_p^{n-d}$. It is denoted with $\underline{s} S_p^{n-d}$ in accordance with the above notation $jS_p^{n-1}$. Note that $S_p^n$ contains $p^{d}$  pairwise disjoint subgraphs $\underline{s} S_p^{n-d}$, $\underline{s}\in [p]_0^{d}$. In particular, $S_p^n$ contains $p^{n-1}$ pairwise disjoint $p$-cliques $\underline{s} S_p^{1}$, $\underline{s}\in [p]_0^{n-1}$. The vertices $i^n$, $i\in [p]_0$, of $S_p^n$ are called {\em extreme vertices} (of $S_p^n$).  The clique in which an extreme vertex lies is called an {\em extreme clique}.

After this preparation we can state the following result which asserts, roughly speaking, that precisely one half of the Sierpi\'nski graphs are EOCD graphs.

\begin{theorem}
\label{thm:eocd-sirpinski}
Let $p\ge 3$ and $n\ge 2$. Then $S_p^n$ is an EOCD graph if and only if $p$ is even.
\end{theorem}

\proof
Suppose that $p$ is odd and that $D$ is an EOD set of $S_p^n$. Observe first that no extreme vertex of $S_p^n$ lies in $D$ because otherwise $D$ would contain two vertices from the same extreme clique, which is not possible. Hence every vertex from $D$ is of degree $p$ and consecutively $|D| = |V(S_p^n)|/p = p^{n-1}$. Since this is at the same time the number of all $p$-cliques of $S_p^n$, it follows that $D$ must have precisely one vertex in each $p$-clique of $S_p^n$. By the same argument as above, a vertex $s$ of $D$ can only be covered by a vertex $t$ of $D$ that lies in a $p$-clique that is neighboring the $p$-clique of $s$. This means that the vertices of $D$ can be partitioned into disjoint pairs. But $p$ is odd and hence $|D| = p^{n-1}$ is odd as well, hence $D$ does not exist.

Assume now that $p$ is even, say $p=2k$, $k\ge 2$. We first recall from~\cite{klavzar-2002} that $S_p^n$ contains an ECD set. In order to prove that $S_p^n$ is an EOCD graph it thus remains to prove that it contains an EOD set. For this sake set
$$D_{2i} = \{\underline{s}(2i)(2i+1):\ \underline{s}\in [p]_0^{n-2}\},\quad 0\le i\le k-1\,,$$
and
$$D_{2i+1} = \{\underline{s}(2i+1)(2i):\ \underline{s}\in [p]_0^{n-2}\},\quad 0\le i\le k-1\,.$$
We claim that
$$D = \bigcup_{i=0}^{2k-1}D_{i}$$
is an EOD set of $S_p^n$. Note first that for any $i\in [k]_0$, $|D_{2i}| = |D_{2i+1}| = p^{n-2}$. Since the sets $D_i$, $i\in [2k]_0$, are clearly pairwise disjoint, it follows that $|D| = 2kp^{n-2} = p^{n-1}$. Let now $\underline{s} S_p^{1}$, $\underline{s}=s_n\ldots s_{2}\in [p]_0^{n-1}$, be an arbitrary $p$-clique of $S_p^n$. If $s_2$ is even, say $s_2=2i$, then $\underline{s}(2i+1)\in D\cap \underline{s} S_p^{1}$, and if $s_2$ is odd, say $s_2=2i+1$, then $\underline{s}(2i)\in D\cap \underline{s} S_p^{1}$. If follows that any $p$-clique contains a vertex of $D$ and consequently it contains exactly one such vertex. Since by the construction of the sets $D_{2i}$ and $D_{2i+1}$ any vertex of $D$ has a neighbor in $D$, we conclude that $D$ is indeed an EOD set of $S_p^n$.
\qed

Combining the construction of the EOC sets in the proof of Theorem~\ref{thm:eocd-sirpinski} with Proposition~\ref{prp:EOD-ECD}(ii) we get:

\begin{corollary}
If $p\ge 4$  is even and $n\ge 2$, then $\gamma_t(S_p^n) = p^{n-1}$.
\end{corollary}

\section*{Acknowledgements}

The work of SK and IP was in part financed by ARRS Slovenia under the grant P1-0297. Most of the research was done while IGY was visiting the University of Maribor, Slovenia, as part of the grant ``Internationalisation---a pillar of development of University of Maribor."


\end{document}